%% file: zhang2.tex
\documentclass[12pt]{amsart}
\usepackage{color,graphicx,wrapfig}
\def\div{\hbox{\rm div }}

\def\R{\mathbb{R}}

\begin{document}

\title[Existence of a degenerate singularity]{Existence of a degenerate singularity in the high activation energy
limit of a reaction-diffusion equation}
\author[G.S. Weiss]{Georg S. Weiss}
\address{Graduate School of Mathematical Sciences,
University of Tokyo, 3-8-1 Komaba, Meguro-ku, Tokyo-to, 153-8914 Japan}
\email{gw@ms.u-tokyo.ac.jp}
\thanks{$2000$ {\it Mathematics Subject Classification.\/} Primary
35R35, Secondary 35J60.}
\thanks{G.S. Weiss has been partially supported by the Grant-in-Aid 18740086 of the Japanese Ministry of Education, Culture, Sports, Science 
and Technology}
\author[G. Zhang]{Guanghui Zhang}
\address{Graduate School of Mathematical Sciences,
University of Tokyo, 3-8-1 Komaba, Meguro-ku, Tokyo-to, 153-8914 Japan}
\email{zhang@ms.u-tokyo.ac.jp}
\date{}
\theoremstyle{plain}\newtheorem*{thma}{Theorem A}
\newtheorem{thm}{Theorem}[section]
\newtheorem{lem}[thm]{Lemma}
\newtheorem{prop}[thm]{Proposition}
\newtheorem{cor}[thm]{Corollary} \theoremstyle{definition}
\newtheorem{defn}[thm]{Definition} \theoremstyle{example}
\newtheorem{example}[thm]{Example}
\theoremstyle{definition}
\newtheorem{remarks}[thm]{Remark}
\numberwithin{equation}{section}
\maketitle

\begin{abstract}
We consider the singular perturbation problem
$$
\Delta u_\epsilon=\beta_\epsilon(u_\epsilon),
$$
where
$\beta_\epsilon(s)=\frac{1}{\epsilon}\beta(\frac{s}{\epsilon})$,
$\beta$ is a Lipschitz continuous function such that $\beta>0$ in
$(0, 1)$, $\beta\equiv 0$ outside $(0, 1)$ and
$\int_0^1\beta(s)\> ds=\frac{1}{2}$.

We construct an example exhibiting a {\em degenerate singularity}
as $\epsilon_k\searrow 0$.
More precisely, there is a sequence of solutions $u_{\epsilon_k}\to u$
as $k\to \infty$, and there exists $x^0\in\partial\{u>0\}$
such that
$$ \frac{u(x^0+r\cdot)}{r} \to 0 \textrm{ as } r\to 0.$$
Known results suggest that this singularity must be {\em unstable},
which makes it hard to capture analytically and numerically.
Our result answers a question raised by Jean-Michel Roquejoffre
at the FBP'08 in Stockholm.
\end{abstract}
\section{Introduction}
In this paper we are concerned with the singular perturbation
problem
\begin{equation}\label{E}
\Delta u_\epsilon=\beta_\epsilon(u_\epsilon) \> \textrm{ in } \Omega,
\end{equation}
where $\Omega$ is a bounded domain in $\R^n$, $\epsilon>0$ and
$\beta_\epsilon(s)=\frac{1}{\epsilon}\beta(\frac{s}{\epsilon})$.
Here $\beta$ is a Lipschitz continuous function such that $\beta>0$ in
$(0, 1)$, $\beta\equiv 0$ outside $(0, 1)$ and
$\int_0^1\beta(s)\> ds=\frac{1}{2}$.

This problem arises in the mathematical analysis of
equidiffusional flames (see \cite{Berestycki}, \cite{Bucmaster}),
in which case $\epsilon$ is proportional to the inverse of the
activation energy.

Formally, as $\epsilon\to 0$, the solutions $u_\epsilon$
converge to a solution $u$ of the free boundary
problem 
\begin{equation}\label{fbp}
\begin{array}{ll}
\Delta u=0  &\hbox{in}\ \Omega\setminus\partial\{u>0\},\\
u=0, \ (u^+_\nu)^2-(u^-_\nu)^2=1 & \hbox{on}\
\partial\{u>0\}.
\end{array}
\end{equation}

On a rigorous level, 
L. Caffarelli established a locally
uniformly Lipschitz estimate for bounded solutions $\{u_\epsilon\}$
(see \cite{Caffarelli}).
In \cite{Lederman}, C. Lederman and
N. Wolanski proved that $u$ is a viscosity solution of 
(\ref{fbp}). They also proved that
$u$ satisfies the free boundary condition in a pointwise sense
at non-degenerate free boundary points at which there is an inward unit normal of
$\{u>0\}$ in the measure theoretic sense. Here ``non-degenerate''
means that
$$r^{-n-1} \int_{B_r(x)} u^+ \ge c>0\textrm{ for } r \le r_0.$$
A related convergence result has been proved for the $p$-Laplace operator in
\cite{DPS}.

From \cite[Lemma 3.4]{AC} we know that all free boundary points
of local minimizers are non-degenerate.
This suggests that if there are degenerate free boundary
points of a limit function $u$, then they must be unstable.

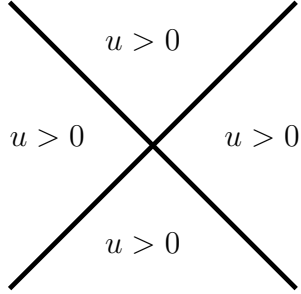
\begin{figure}
\begin{center}
\input{zhang1.pstex_t}
\end{center}
\caption{Cross-like Singularity}\label{zhang1}
\end{figure}

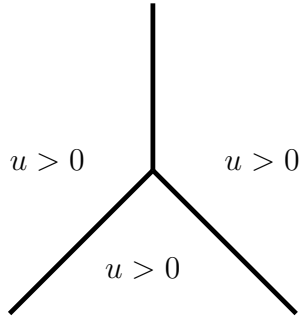
\begin{figure}
\begin{center}
\input{zhang2.pstex_t}
\end{center}
\caption{Triple Junction Singularity}\label{zhang2}
\end{figure}

It has been known that cross-shaped free boundaries occur
for domain variation solutions (see the Introduction of
\cite{Weiss2}),
for example
$$ u(x_1,x_2)=|x_1^2-x_2^2| \textrm{ (cf. Figure \ref{zhang1}).}$$
\newline
John Andersson suggested a triple junction example
(see Figure \ref{zhang2})
where the solution is homogeneous of degree
$3/2$; domain variation solutions with
a degree of homogeneity $\in (1, 3/2)$ are not possible.
\\  
However it is not so obvious whether these examples arise as
limits of the reaction-diffusion equation (\ref{E}).
This is the question
Jean-Michel Roquejoffre posed
at the FBP'08 in Stockholm.
The main result of the present paper gives a partial answer
to that question:

\begin{thma}
There exist uniformly bounded solutions $u_{\epsilon_k}:B_1(0)\subset \R^2\to \R$
such that $\epsilon_k\to 0, u_{\epsilon_k}\to u$ as $k\to\infty$, and
there is a degenerate free boundary point of $u$ in $B_1(0)$.
More precisely,
there is $x^0\in B_1(0)\cap \partial\{u>0\}$ such that
$$\lim_{r\to
0}\frac{u(x^0+rx)}{r}=0 \textrm{ for every }x\in \mathbb{R}^2.$$
\end{thma}

\begin{figure}
\begin{center}
\input{zhang3.pstex_t}
\end{center}
\caption{Imaginable $\epsilon$-solution}\label{zhang3}
\end{figure}
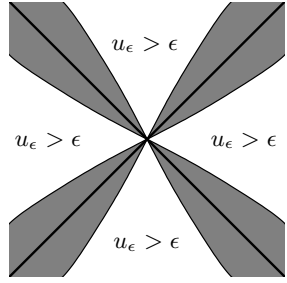

Remark \ref{fin} shows that we may construct degenerate points
with an arbitrarily high number of symmetry lines.
See Figure \ref{zhang3} for an example with 4 symmetry lines.
\newline
Section \ref{3} also suggests a {\em numerical} approach 
for the construction of this unstable solution
which would be otherwise hard to capture.
\section{Notation}
\par
\indent
Throughout this article $\mathbb{R}^n$ will be equipped with the Euclidean
inner product $x\cdot y$ and $B_r(x^0)$ will denote the open ball of center
$x^0$ with  radius $r$. When not specified, $x^0$ is assumed to be $0$. 
We will use the polar coordinates
$(r, \theta)$ in $\mathbb{R}^2$.  We use 
the $n$-dimensional Lebesgue measure $\mathcal{L}^n$ and the $m-$dimensional 
Hausdorff measure $\mathcal{H}^m$. We shall often use abbreviations for inverse images like
$\{u>0\}:=\{x\in\Omega : u(x)>0\}$.

\section{Construction of $\epsilon-$solutions}\label{3}
We are going to use a method introduced in \cite{AW}.
\par
Let $K=\{(r\cos\theta, r\sin\theta):0<r<1, 0<\theta<\pi/m\}$,
where $m$ is an integer $\ge 3$. For $g\in
C^{\alpha_0}(\partial B_1\cap\partial K)$, $\alpha_0\ge \alpha >0$, $C_g^{\alpha}(\bar{K})$ will
denote the subspace of $C^{\alpha}(\bar{K})$ consisting of all
functions with boundary values $g$ on $\partial B_1\cap\partial K$.
We consider the following problem:

\begin{equation}\label{E2}
\begin{array}{l l}
\Delta u(x)=\beta_\epsilon(u(x)-u(0)+\epsilon)\quad &\mbox{in}\  K, \\
u=g   &\mbox{on}\  \partial K\cap\partial B_1, \\
\frac{\partial u}{\partial\nu}=0  &\mbox{on}\  \partial
K\setminus\partial B_1.
\end{array}
\end{equation}

\par
Let $T=T_{\epsilon, g}:C_g^\alpha(\bar{K})\to C_g^\alpha(\bar{K})$ be
the operator defined by
\begin{equation*}
\begin{array}{l l}
\Delta T(u)=\beta_\epsilon(u(x)-u(0)+\epsilon)\quad &\mbox{in}\ K, \\
T(u)=g  &\mbox{on}\  \partial K\cap\partial B_1, \\
\frac{\partial T(u)}{\partial\nu}=0  & \mbox{on}\  \partial
K\setminus\partial B_1.
\end{array}
\end{equation*}

\par
Using Schauder's fixed point theorem,  we can prove the existence
of a solution of equation (\ref{E2}):
\begin{prop}
For each $g\in C^{\alpha_0}(\partial B_1\cap\partial K)$, T has a fixed
point $u_\epsilon\in C^\alpha_g(\bar{K})$ for some $\alpha>0$
depending only on $\alpha_0$.
\begin{proof}
Let $v\in C^\alpha_g(\bar{K})$, $F(x)=\beta_\epsilon(v(x)-v(0)+\epsilon)$. Then there is a $W^{1,2}(K)$-solution $u$ of the problem:
\begin{equation*}
\begin{array}{l l}
\Delta u=F\quad  & \mbox{in}\  K, \\
u=g  & \mbox{on}\ \partial K \cap\partial B_1, \\
\frac{\partial u}{\partial\nu}=0  &  \mbox{on}\  \partial
K\setminus\partial B_1.
\end{array}
\end{equation*}
We extend $u$ to a function $\tilde{u}:B_1\setminus \{ 0\}\to \R$ by even reflection:
\begin{displaymath}
\begin{array}{l l}
\tilde{u}(r, \pi/m+\theta)=u(r, \pi/m-\theta), \quad  &\mbox{if}\
0\leq\theta\leq \pi/m,\\
\tilde{u}(r, j\pi/m+\theta)=\tilde{u}(r, j\pi/m-\theta), \quad  &\mbox{if}\
0\leq\theta\leq \pi/m,\\
\end{array}
\end{displaymath}
for $j=1, 2, \dots, 2m-1$. It follows that $\tilde{u}$ is a solution of
\begin{equation*}
\begin{array}{l l}
\Delta\tilde{u}=\tilde{F}\quad  & \mbox{in}\  B_1\setminus\{0\}, \\
\tilde{u}=\tilde{g}  & \mbox{on}\  \partial B_1,
\end{array}
\end{equation*}
where $\tilde{F}$ and $\tilde{g}$ are defined by even reflection.
As the origin is a set of vanishing 2-capacity (see \cite{Frehse}),
$\tilde u$ is a weak solution of $\Delta \tilde u=\tilde{F}$ in $B_1$.
\par 
Applying the regularity theory of elliptic equations (see 
for example \cite[Corollary 9.29]{Gilbarg}), we see that
$\tilde{u}\in C^\alpha(\bar{B}_1)$ for small $\alpha$, and that $||\tilde{u}||_{C^{0,\alpha}(B_1)}\leq C$,  where $\alpha\in (0,\alpha_0)$ and $C$
are constants depending on $\alpha_0$, $||g||_{C^{0, \alpha_0}}$ and 
$||F||_{L^\infty(K)}$. Hence $T$ is a continuous compact linear
operator from $C_g^\alpha(\bar{K})$ into itself for small
$\alpha>0$, and
$$
||T_{\epsilon, g}(w)||_{C^\alpha(\bar{K})}\leq C
$$
for every $\sigma\in [0,1]$ and every solution of the equation
\begin{equation*}
\begin{array}{l l}
\Delta w=\sigma\beta_\epsilon(w(x)-w(0)+\epsilon)\quad  &
\mbox{in}\ K,\\
w=g  &\mbox{on}\  \partial K \cap\partial B_1,\\
\frac{\partial w}{\partial\nu}=0  &\mbox{on}\  \partial
K\setminus\partial B_1;
\end{array}
\end{equation*}
here $C$ is a constant depending only on $g$, $\beta$ and $\epsilon$.
\par
From Schauder's fixed point theorem (see \cite[Chapter 11]{Gilbarg})
we infer that $T_{\epsilon, g}$ has a fixed point $u_\epsilon\in C^\alpha_g(\bar{K})$.
\end{proof}
\end{prop}
\section{Convergence}
From now on all of our functions are defined in $B_1(0)$ by the
above reflection. We choose one non-constant $g\in C^2(\bar{B})$.
\par
Let $v_\epsilon(x)=u_\epsilon(x)-u_\epsilon(0)+\epsilon$. Then $v_\epsilon$ is
a solution of
\begin{equation*}
\begin{array}{l l}
\Delta v_\epsilon(x)=\beta_\epsilon(v_\epsilon(x))\quad & \hbox{in}\ B_1,\\
v_\epsilon(0)=\epsilon,\\
v_\epsilon=g-u_\epsilon(0)+\epsilon &  \hbox{on}\  \partial B_1.\\
\end{array}
\end{equation*}
By the maximum principle, $|u_\epsilon(0)|\leq||g||_{L^\infty(B_1)}$,
therefore $||v_\epsilon||_{L^\infty(\bar{B}_1)}\leq C$, where $C$ depends only on
$||g||_{L^\infty(B_1)}$.
\par
First we state a result from \cite{Caffarelli} proving a uniform Lipschitz
estimate.

\begin{prop}\label{uniformly Lipschitz}
Let $w_\epsilon$ be a family of solutions to $\Delta
w_\epsilon=\beta_\epsilon(w_\epsilon)$ in a domain $\Omega\subset
\mathbb{R}^n$ such that $||w_\epsilon||_{L^\infty(\Omega)}\leq C$ for some
$C<\infty$. Let $K\subset\Omega$ be a compact set and let $\tau>0$ be such
that $B_\tau(x^0)\subset\Omega$ for every $x^0\in K$. Then there
exists a constant $L=L(\tau, C)$, such that
$$
|\nabla w_\epsilon(x)|\leq L \ \ \ \mathrm{for}\  x\in K.
$$
\end{prop}
If the boundary data is smooth enough, we have the following gradient
estimate near the boundary (see \cite{Gurevich}).
\begin{prop}\label{boundary Lipschitz}Let $\Omega\subset\mathbb{R}^n$ be a bounded domain with smooth boundary 
 and let $w_\epsilon\in C^2(\Omega)\cap C^0(\bar{\Omega})$ be a solution of
\begin{equation*}
\begin{array}{l l}
\Delta w_\epsilon(x)=\beta_\epsilon(w_\epsilon(x))\quad &\mbox{in}\ \Omega,\\
w_\epsilon=f &\mbox{on}\  \partial \Omega,
\end{array}
\end{equation*}
where $f\in C^2(\bar{\Omega})$, $||f||_{C^2(\bar{\Omega})}<C_1$. Then there
exist $\epsilon_0>0, C<\infty$ such that
$$
|\nabla w_\epsilon(x)|\leq C(1+|\log d(x, \partial\Omega)|)
$$
for $\epsilon \in (0,\epsilon_0)$, $x\in\Omega$, where $d(x, \partial\Omega)$ is the 
distance from $x$ to the boundary of $\Omega$.
\end{prop}
Now we are ready to prove the convergence of $\{v_{\epsilon_i}\}$.
\begin{prop}\label{conv}
There exist a sequence $\{v_{\epsilon_i}\}$ and a locally
Lipschitz continuous function $v_0$ such that\\
1) $v_{\epsilon_i}\rightarrow v_0$ uniformly on $\bar{B}_1$, \\
2) $v_{\epsilon_i}\rightarrow v_0$ in $W^{1, 2}_\mathrm{loc}(B_1)$, \\
3) $v_0$ is harmonic in $B_1\setminus\partial\{v_0>0\}$, \\
4) $\Delta v_{\epsilon_i}\rightarrow \mu$ as measures on $B_1$; here
$\mu$ is a locally finite non-negative measure supported on the free boundary
$\partial\{v_0 >0\}$.Therefore $$ \Delta v_0=\mu\qquad \hbox{in}\
B_1, $$that is
$$\int_{B_1}\nabla v_0\cdot \nabla\phi\ dx+\int_{B_1}\phi\ d\mu=0$$
for every $\phi\in C_0^1(B_1)$.
\\
5) $v_0(r, j\pi/m+\theta)=v_0(r, j\pi/m-\theta)$ for
$0\leq\theta\leq \pi/m$ and $1 \le j \le 2m$.
\begin{proof}
1) By Proposition \ref{boundary Lipschitz}, $\{v_\epsilon\}$ is
equicontinuous. Since $\{v_\epsilon\}$ is uniformly bounded in
$C^0(\bar{B}_1)$, by the Arzela-Ascoli theorem there is a sequence
$\{v_{\epsilon_i}\}$ that converges to some
$v_0$ in $C^0(\bar{B}_1)$. By Proposition \ref{uniformly Lipschitz},
it is easy to see that $v_0$ is locally Lipschitz continuous. \\

2) For $0\leq\eta\in C_0^\infty(B_1)$, we obtain from the local Lipschitz bound that
\begin{equation}\label{ineq1}
\int_{B_1}\eta|\nabla v_0|^2\ dx\leq\liminf_{\epsilon_i\to
0}\int_{B_1}\eta|\nabla v_{\epsilon_i}|^2\ dx. 
\end{equation}
If we prove the converse inequality, the result will
follow. Multiplying the $\epsilon_i$-equation (\ref{E}) by $\eta
v_{\epsilon_i}$, we obtain that

\begin{equation*}
\int_{B_1}v_{\epsilon_i}\nabla\eta\cdot\nabla
v_{\epsilon_i}\ dx+\int_{B_1}\eta|\nabla
v_{\epsilon_i}|^2\ dx=-\int_{B_1}\eta
v_{\epsilon_i}\beta_{\epsilon_i}(v_{\epsilon_i})\ dx\leq 0.
\end{equation*}
Letting $\epsilon_i\to 0$, we get
\begin{equation}\label{ineq2}
\limsup_{\epsilon_i\to0}\int_{B_1}\eta|\nabla
v_{\epsilon_i}|^2\ dx\leq-\int_{B_1} v_0\nabla\eta\cdot\nabla v_0\ dx.
\end{equation}
For each $\delta>0$, we define
\begin{equation*}
v_\delta(x)=\begin{cases} v_0(x)-\delta, & \hbox{if}\  v_0(x)>\delta, \\
v_0(x)+\delta, &\hbox{if}\  v_0(x)<-\delta, \\
0 & \hbox{otherwise}.
\end{cases}
\end{equation*}
Multiplying the $\epsilon_i$-equation by $\eta v_\delta$ and
integrating, we have
\begin{equation*}
\int_{B_1}v_\delta\nabla\eta\cdot\nabla
v_{\epsilon_i}\ dx+\int_{B_1}\eta\nabla v_\delta\cdot\nabla
v_{\epsilon_i}\ dx=0.
\end{equation*}
Letting $\epsilon_i\to 0$, it follows that
\begin{equation*}
\int_{B_1}v_\delta\nabla\eta\cdot\nabla v_0\ dx+\int_{B_1}\eta|\nabla
v_\delta|^2\ dx=0,
\end{equation*}
and letting $\delta\to 0$, we get
\begin{equation*}
\int_{B_1}v_0\nabla\eta\cdot\nabla v_0\ dx+\int_{B_1}\eta|\nabla
v_0|^2\ dx=0.
\end{equation*}
Using \eqref{ineq2}, we obtain
\begin{equation}\label{ineq3}
\limsup_{\epsilon_i\to0}\int_{B_1}\eta|\nabla
v_{\epsilon_i}|^2\ dx\leq\int_{B_1} \eta|\nabla v_0|^2\ dx.
\end{equation}
Combining (\ref{ineq1}) and (\ref{ineq3}), we obtain
\begin{equation}
\lim_{\epsilon_i\to0}\int_{B_1}\eta |\nabla
v_{\epsilon_i}|^2\ dx=\int_{B_1}\eta|\nabla v_0|^2\ dx,
\end{equation}
which implies that
$$\eta^{1/2}\nabla
v_{\epsilon_i}\to \eta^{1/2}\nabla
v_0\ \ \hbox{in}\ L^2(B_1). $$
For $\eta$ such that $\eta=1$ in $D\subset\subset\Omega$, it follows that 
$$\nabla v_{\epsilon_i}\to \nabla
v_0\ \ \hbox{in}\ L^2(D).$$
Therefore $v_\epsilon\to v_0$ strongly in $W^{1, 2}(D)$.
\par
By a diagonal sequence argument, assertion 2) of the proposition
follows. \\

3) Since $v_0$ is continuous,  the sets $\{v_0>0\}$ and
$\{v_0<0\}$ are open. Let $x^0\in\{v_0>0\}$. From the fact that
$v_{\epsilon_i}\rightarrow v_0$ uniformly on $B_1$, there exist
$r>0$ and $N\in \mathbb{N}$ such that $v_{\epsilon_i}(x)\geq
v_0(x^0)/2>0$ for $x\in B_r(x^0)$, $i\geq N$. Thus $v_{\epsilon_i}$
is harmonic in $B_r(x^0)$ for $i\geq N$ and the same fact holds for
$v_0$. In the same way we prove that $v_0$ is superharmonic
in $\{v_0\leq 0\}^0$. On the other hand, $v_{\epsilon_i}$ is
subharmonic in $B_1$, so that $v_0$ is also subharmonic.
Therefore $v_0$ is harmonic in $\{v_0\leq 0\}^0$ and
3) follows. \\

4) Multiplying equation $\Delta
v_{\epsilon_i}=\beta_{\epsilon_i}(v_{\epsilon_i})$ by $\phi\in
C_0^\infty(B_1)$ and integrating by parts, we get
$$\int_{B_1}\nabla v_{\epsilon_i}\cdot\nabla\phi
\ dx+\int_{B_1}\beta_{\epsilon_i}(v_{\epsilon_i})\phi\ dx=0. $$ We
know that $v_{\epsilon_i}\rightarrow v_0$ in $W^{1,
2}_\mathrm{loc}(B_1)$, hence
\begin{equation}\label{eq1}
\int_{B_1}\nabla v_{\epsilon_i}\cdot\nabla\phi\
dx\rightarrow\int_{B_1}\nabla v_0\cdot\nabla\phi\ dx. 
\end{equation}
On the other hand, 
$$\int_{B_1}\beta_{\epsilon_i}(v_{\epsilon_i})\phi\ dx\leq C.$$
This $L^1$-bound implies that there exists a locally finite non-negative measure
$\mu$ such that for a subsequence which we still call
$\{v_{\epsilon_i}\}$, $\beta_{\epsilon_i}(v_{\epsilon_i})\rightarrow\mu$ as measures in
$B_1$. Passing to the limit in (\ref{eq1}), we get
$$\int_{B_1}\nabla v_0\cdot\nabla\phi\ dx+\int_{B_1}\phi\ d\mu=0, $$
which implies $$ \Delta v_0=\mu\qquad \hbox{in}\ B_1. $$ Since we know that
$v_0$ is harmonic in $B_1\setminus\partial\{v_0>0\}$, we conclude
that $\hbox{supp}\ \mu\subset\partial\{v_0>0\}$.\\

5) follows from the fact that $v_{\epsilon_i}$ has been defined via even
reflection. This property is preserved by the uniform limit $v_0$.
\end{proof}
\end{prop}
\par
Let
$\mathcal{B}_{\epsilon}(z)=\int_0^z\beta_{\epsilon}(s)\> ds$ and $\chi_\epsilon(x)=2\mathcal{B}_\epsilon(v_\epsilon(x))$. It follows that $0\leq\chi_{\epsilon}(x)\leq
1$, 
and the uniformly Lipschitz estimate for $v_\epsilon$ implies the relative compactness of
$\chi_{\epsilon}$ in $L_{\mathrm{loc}}^1(B_1)$:
\begin{prop}
$\{\chi_{\epsilon}\}$ is precompact in $L^1(D)$ for each $D\subset\subset B_1$.
\begin{proof}
By Proposition \ref{uniformly Lipschitz}, there exists an $L>0$ such
that $||\nabla v_\epsilon||_{L^\infty(D)}\leq L$, hence
$$
\int_D|\nabla\chi_{\epsilon}|\ dx=2\int_D|\beta_\epsilon(v_\epsilon)\nabla
v_\epsilon|\ dx\leq
2L\int_D\beta_\epsilon(v_\epsilon)\ dx\leq C.
$$
Therefore $\{\chi_{\epsilon}\}$ is bounded in $W^{1, 1}(D)$.
By the Sobolev embedding theorem $\{\chi_{\epsilon}\}$ is
precompact in $L^1(D)$.
\end{proof}
\end{prop}
Thus we may assume that
$\chi_{\epsilon_i}(x)\rightarrow\chi_0(x)$ locally in $L^1$.
Similar to \cite[Lemma 4.1]{Weiss2}, we have
\begin{prop}
$\chi_0(x)\in\{0, 1\}$ for a.e. $x\in B_1$.
\begin{proof}
For $\delta\in(0, 1/4)$, $\mathcal{B}(s)=2\int_0^s\beta(t)\ dt$ and
$A_\delta=\{x\in D\subset\subset B_1:\chi_0\in[2\delta, 1-2\delta]\}$ we have
\begin{align*}
\mathcal{L}^2(A_\delta)&\leq\mathcal{L}^2(A_\delta\cap\{\chi_{\epsilon_m}\not\in[\delta, 1-\delta]\})
+\mathcal{L}^2(A_\delta\cap\{\mathcal{B}(\frac{v_{\epsilon_m}}{\epsilon_m})\in[\delta, 1-\delta]\})\\
&\leq\mathcal{L}^2(D\cap \{|\chi_{\epsilon_m}-\chi_0|>\delta\})+
\mathcal{L}^2(D\cap \{\frac{v_{\epsilon_m}}{\epsilon_m}\in[\mathcal{B}^{-1}(\delta), \mathcal{B}^{-1}(1-\delta)]\})\\
&\leq\mathcal{L}^2(D\cap \{|\chi_{\epsilon_m}-\chi_0|>\delta\})+C_1(\delta, \beta)\epsilon_m\int_D\beta_{\epsilon_m}(v_{\epsilon_m})\ dx\\
&\leq\mathcal{L}^2(D\cap \{|\chi_{\epsilon_m}-\chi_0|>\delta\})+C_2(\delta, \beta)\epsilon_m\rightarrow 0
\end{align*}
as $m\to\infty$.
\end{proof}
\end{prop}
\section{Monotonicity formula}
In this section we state a monotonicity formula proved in \cite{Weiss1}, 
which is a key tool in constructing a degenerate singularity.

\begin{prop}\label{monotonicity formula}
Let $x^0\in B_1$, and
$$\Phi_{\epsilon}(r)=r^{-2}\int_{B_r(x^0)}|\nabla{v_{\epsilon}}|^2+r^{-2}\int_{B_r(x^0)}\chi_{\epsilon}-r^{-3}\int_{\partial{B_r(x^0)}}v_{\epsilon}^2\ d\mathcal{H}^1. $$
Then $\Phi_\epsilon$ satisfies the monotonicity formula
$$\Phi_{\epsilon}(\sigma)-\Phi_{\epsilon}(\rho)\geq\int_\rho^\sigma r^{-2}\int_{\partial B_r(x^0)}2(\nabla v_\epsilon\cdot\nu-\frac{v_\epsilon(x)}{r})^2\ d\mathcal{H}^1dr\ \ for\
0<\rho<\sigma.
$$
Letting $\epsilon\to 0$,
$$\Phi(r)=r^{-2}\int_{B_r(x^0)}|\nabla{v_0}|^2+r^{-2}\int_{B_r(x^0)}\chi_0-r^{-3}\int_{\partial{B_r(x^0)}}{v_0}^2\ d\mathcal{H}^1$$
satisfies the monotonicity formula
$$\Phi(\sigma)-\Phi(\rho)\geq\int_\rho^\sigma r^{-2}\int_{\partial B_r(x^0)}2(\nabla v_0\cdot\nu-\frac{v_0(x)}{r})^2\ d\mathcal{H}^1dr\ \ for\
0<\rho<\sigma. $$
\end{prop}

\section{Existence of a Degenerate Free Boundary Point}
In this section we are going to prove the existence of a degenerate singular free 
boundary point for $v_0$. We start with some lemmas.
\begin{lem}\label{domain variation}
Let $w$ be a solution of
$$
\Delta w=\beta_\epsilon(w)\textrm{ in }
\Omega\subset\mathbb{R}^n\textrm{ and let } \chi_\epsilon(x) = 2\mathcal{B}_\epsilon(w(x)).
$$
Then 
$$
\int_\Omega|\nabla w|^2 \div\phi-2\nabla w D\phi\nabla
w+\chi_\epsilon\div\phi\ dx=0 
$$ 
for $\phi\in C_0^1(\Omega;\mathbb{R}^n)$.
\begin{proof}
Integrating by parts, we get
\begin{equation*}
\begin{split}
&\int_\Omega|\nabla w|^2 \div\phi-2\nabla w D\phi\nabla
w+\chi_\epsilon\div\phi\ dx\\
&=\int_\Omega\sum_{i=1}^n\left(\frac{\partial w}{\partial
x_i}\right)^2\sum_{j=1}^n\frac{\partial \phi_j}{\partial
x_j}-2\sum_{i, j=1}^n\frac{\partial w}{\partial x_i}\frac{\partial
\phi_j}{\partial x_i}\frac{\partial w}{\partial
x_j}+\chi_\epsilon\div\phi\ dx\\
&=-2\int_\Omega\sum_{i, j=1}^n\phi_j\frac{\partial w}{\partial
x_i}\frac{\partial^2 w}{\partial x_i\partial
x_j}-\sum_{i, j=1}^n\phi_j\frac{\partial}{\partial
x_i}\left(\frac{\partial w}{\partial x_i}\frac{\partial w}{\partial
x_j}\right)+\beta_\epsilon(w)\nabla w\cdot\phi\ dx\\
&=2\int_\Omega\sum_{i, j=1}^n\phi_j\frac{\partial w}{\partial
x_j}\frac{\partial^2w}{\partial x_i^2}-\beta_\epsilon(w)\nabla w\cdot\phi\ dx\\
&=2\int_\Omega\nabla w\cdot\phi(\Delta w-\beta_\epsilon(w))\ dx\\
&=0. 
\end{split}
\end{equation*}
\end{proof}
\end{lem}

\begin{lem}\label{chi_limit}
If $v_0\equiv 0$ in $B_r$ for some $r>0$, then $\chi_0\equiv 1$ in
$B_r$.
\begin{proof}
By Lemma \ref{domain variation} we have
$$
\int_{B_1}|\nabla v_{\epsilon_i}|^2 \div\phi-2\nabla v_{\epsilon_i}D\phi\nabla
v_{\epsilon_i}+\chi_{\epsilon_i}\div\phi\ dx=0
$$ for
$\phi\in C_0^\infty(B_r; \mathbb{R}^2)$. Letting $\epsilon_i\to 0$, we obtain from 
the strong convergence of $v_{\epsilon_i}$ that 
\begin{align*}
&0=\int_{B_1}|\nabla v_0|^2 \div\phi-2\nabla v_0D\phi\nabla
v_0+\chi_0 \div\phi\ dx\\
&=\int_{B_r}\chi_0 \div\phi\ dx,
\end{align*}
hence $\chi_0\equiv\hbox{const}$ in $B_r$. We know that
$\chi_0(x)\in\{0, 1\}$ for a.e. $x\in B_1$, therefore
$\chi_0(x)\equiv0$ or $\chi_0(x)\equiv1$ in $B_r$.
\par
For $0<c\le 1/4$, let $P_i=\{x\in B_r:c\epsilon_i\leq v_{\epsilon_i}(x)\leq
(1-c)\epsilon_i\}$ and $Q_i=\{x\in B_r:v_{\epsilon_i}(x)>c\epsilon_i\}$.
Supposing towards a contradiction that $\chi_0\equiv 0$ in
$B_r$, and recalling that $\chi_{\epsilon_i}\rightarrow\chi_0$
in $L^1(B_r)$ we have
$$
\mathcal{L}^2(Q_i)\leq\mathcal{L}^2(\{x:\chi_{\epsilon_i}>2\mathcal{B}(c)\})\rightarrow 0
$$
as $i\to\infty$.
\par
Let $\delta\in (0,r/4)$, $I_{\epsilon_i}=\{t :\delta\leq t\leq
r,\min\limits_{x\in\partial B_t}v_{\epsilon_i}(x)\geq \epsilon_i/2\}$ and
$I^c_{\epsilon_i}=[\delta, r]\setminus I_{\epsilon_i}$. Then we have
$\mathcal{L}^1(I_{\epsilon_i})\leq\mathcal{L}^2(Q_i)/\delta\to 0$ as
$i\to\infty$, hence there exists an $\epsilon_0>0$ such that
$\mathcal{L}^1(I_{\epsilon_i}^c)>r/2$ for $\epsilon_i<\epsilon_0$. As $v_{\epsilon_i}$ is subharmonic in $B_1$,
$\epsilon_i= v_{\epsilon_i}(0)\le\max\limits_{\partial B_t}v_{\epsilon_i}$ for every $t\in[ 0, r)$. It follows that 
for each $t\in I_{\epsilon_i}^c$, there exists a point $x^i_t\in\partial
B_t$ such that $v_{\epsilon_i}(x_t^i)=\epsilon_i/2$ (cf. Figure \ref{zhang4}). 
\begin{figure}
\begin{center}
\input{zhang4.pstex_t}
\end{center}
\caption{The set $P_i$}\label{zhang4}
\end{figure}
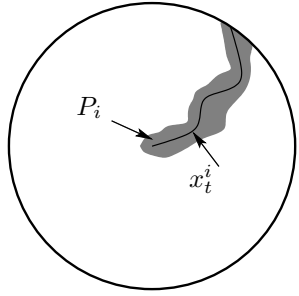
By the
uniformly Lipschitz continuity of $v_{\epsilon_i}$ there exists a
$d>0$ independent of $\epsilon$, such that $B_{ d\epsilon_i}(x^i_t)\subset P_i$. Therefore
$$\mathcal{L}^2(P_i)\geq\mathcal{L}^2(\bigcup_{x_t^i\in I^c_{\epsilon_i}} B_{d\epsilon_i}(x^i_t))\geq rd\epsilon_i/4. $$
On the other hand, since $\Delta v_{\epsilon_i}\rightarrow 0$ in $B_r$,
$$
\mathcal{L}^2(P_i)/\epsilon_i\leq
C\int_{B_r}\beta_{\epsilon_i}(v_{\epsilon_i})\rightarrow0. $$ Thus
we get a contradiction and therefore $\chi_0\equiv 1$ in
$B_r$.
\end{proof}
\end{lem}

\begin{lem}\label{blow-up}
Let $v_0(x^0)=0$ for some $x^0\in B_1$, $0<r_k\rightarrow 0$ as
$k\to\infty$ and $v_{0k}(x)=\frac{v_0(x^0+r_kx)}{r_k}$. Then there
exists a blow-up limit $v_{00}$ such that for a subsequence the following holds:\\
1) $v_{0k}\to v_{00}$ in $C_\mathrm{loc}^{0, \alpha}(\mathbb{R}^2)$ for
every $0<\alpha<1$, $\nabla v_{0k}\to\nabla v_{00}$ weakly* in
$L^{\infty}_{\mathrm{loc}}(\mathbb{R}^2)$, \\
2) $v_{0k}\to v_{00}$ in $W^{1, 2}_{\mathrm{loc}}$, \\
3) $v_{00}$ is homogeneous of degree $1$.
\begin{proof}
1) Since $\nabla v_{0k}(x)=\nabla v_0(x_0+r_kx)$, $\{v_{0k}\}$ is locally
uniformly Lipschitz continuous in $\mathbb{R}^2$ and 1)
follows.\\

\par
2) Let $D\subset\subset\mathbb{R}^2$. Then there exists a
$C<\infty$ such that $|\nabla v_{0k}(x)|<C$ in $D$,  hence $v_{0k}$ is bounded
$W^{1, 2}(D)$. It follows that there is a subsequence, which we still
call $v_{0k}$, such that $v_{0k}\to v_{00}$ weakly in $W^{1, 2}(D)$ and
$$
||v_{00}||_{W^{1, 2}(D)}\leq\liminf_{k\to\infty}||v_{0k}||_{W^{1, 2}(D)}.
$$

\par
Since $\Delta v_0=\mu$ and
$supp\ \mu\subset\partial\{v_0>0\}$, we get
$$
v_{0k}(x)\Delta v_{0k}(x)=0.$$
Let $\eta\in C_0^{\infty}(\mathbb{R}^2)$. Then $\int_{\mathbb{R}^2}\eta|\nabla
v_{0k}|^2=-\int_{\mathbb{R}^2}v_{0k}\nabla\eta\cdot\nabla v_{0k}$ and
$$
\limsup_{k\to0}\int_{\mathbb{R}^2}\eta|\nabla v_{0k}|^2\leq-\int_{\mathbb{R}^2}v_{00}\nabla\eta\cdot\nabla
v_{00}. 
$$ 
For each $\delta>0$, we define
\begin{equation*}
v_\delta(x)=\begin{cases} v_{00}(x)-\delta, & \mbox{if}\  v_{00}(x)>\delta, \\
v_{00}(x)+\delta, &\mbox{if}\  v_{00}(x)<-\delta, \\ 0 & \mbox{otherwise}.
\end{cases}
\end{equation*}
Multiplying the equation for $v_{0k}$ by $\eta
v_{\delta}$, we get
$$
\eta(x)v_{\delta}(x)\Delta v_{0k}(x)=0,
$$ 
if $k$ is large enough. Hence 
$$
\int_{\mathbb{R}^2}\eta\nabla v_{\delta}\cdot\nabla v_{0k}=-\int_{\mathbb{R}^2}v_\delta\nabla\eta\cdot\nabla v_{0k}. 
$$
Letting $k\to\infty$, it follows that
$$
\int_{\mathbb{R}^2}\eta\nabla v_{\delta}\cdot\nabla v_{00}=-\int_{\mathbb{R}^2}v_\delta\nabla\eta\cdot\nabla
v_{00},  
$$
and letting $\delta\to 0$, we get
$$
\int_{\mathbb{R}^2}\eta|\nabla v_{00}|^2=-\int_{\mathbb{R}^2}v_{00}\nabla\eta\cdot\nabla v_{00}. 
$$
Hence 
$$
\limsup_{k\to0}\int_{\mathbb{R}^2}\eta|\nabla
v_{0k}|^2\leq\int_{\mathbb{R}^2}\eta|\nabla v_{00}|^2
$$ 
for every $\eta\in C_0^\infty(\mathbb{R}^2)$. It follows that $v_{0k}\to
v_{00}$ in $W^{1, 2}(D)$.\\
\par 
3) For $0<R<S<\infty$,
$$
\Phi(Rr_k)=R^{-2}\int_{B_R}|\nabla
v_{0k}|^2+\chi_0(r_kx)\ dx-R^{-3}\int_{\partial{B_R}}{v_{0k}}^2\ d
\mathcal{H}^1.
$$
Since $v_0(x^0)=0$ and $v_0$ is Lipschitz continuous, $\Phi(r)$ is bounded. 
Consequently we obtain from the monotonicity formula Proposition \ref{monotonicity formula} that
\begin{align*}
0\leftarrow\Phi(Sr_k)-\Phi(Rr_k)&=\int_R^S2r^{-2}\int_{\partial{B_r}}(\nabla
v_{0k}\cdot\nu-\frac{v_{0k}}{r})^2d\mathcal{H}^1dr\\
&=\int_{B_S\setminus B_R}2|x|^{-4}(\nabla v_{0k}(x)\cdot x-v_{0k}(x))^2dx.
\end{align*}
Letting $k\to\infty$, we obtain that $\nabla v_{00}(x)\cdot x=v_{00}(x)$ a.e.
in $\mathbb{R}^2$, hence $v_{00}$ is homogeneous of degree $1$.
\end{proof}
\end{lem}

Now we are ready to prove Theorem A of the Introduction.
\begin{thm}\label{main}
There is a free boundary point $x^0\in\partial\{v_0>0\}$ such that
$v_0$ is degenerate at $x^0$. More precisely, $$\lim_{r\to
0}\frac{v_0(x^0+rx)}{r}=0$$ for every $x\in \mathbb{R}^2$.
\end{thm}
\begin{proof}
Case 1: $0\not\in\partial\{v_0>0\}$.
\par
There exists an $r>0$ such that $v_0(x)\leq0$ in $B_r$. Since
$v_0(0)=0$, the subharmonicity of $v_0$ implies that $v_0\equiv0$
in $B_r$. Suppose that $r$ is the largest number such that
$v_0\equiv0$ in $B_r$. Noticing that $r<1$ as the boundary
values of $v_0$ are not constant. It follows that there is a point $x^0\in\partial B_r$
such that $x^0\in\partial\{v_0>0\}$ (cf. Figure \ref{zhang5}).
\begin{figure}
\begin{center}
\input{zhang5.pstex_t}
\end{center}
\caption{Touching the Free Boundary}\label{zhang5}
\end{figure}
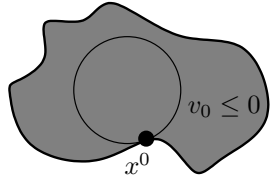
\par
By a translation and a rotation we may assume that $x^0=0$ and $v_0\equiv0$ in
$B_r(-r, 0)$. Let $v_{00}$ be a blow-up limit of $v_0$ at $x^0$, i.e. the limit of $\frac{v_0(x^0+r_kx)}{r_k}$ as $k\to\infty$. 
We are going to show that $v_{00}\equiv0$.
\par
Suppose towards a contradiction that $v_{00}\not\equiv0$ and let
$S=\{ v_{00}>0\}$. Then $v_{00}$ is harmonic in $S$. Moreover, we know from Lemma
\ref{blow-up}  that $v_{00}$ is homogeneous of degree $1$,
so that, solving the resulting ODE for $v_{00}$, $v_{00}=cx\cdot\nu$ in $S$ for some $c\in\mathbb{R}$ and
$\nu\in\partial B_1$. Noticing that $v_{00}\equiv0$ in the left half
plane, we conclude that $v_{00}=c\max(0, x\cdot\nu)$, where $c>0$
and $\nu=(1, 0)$.
\par
Let $\chi_r(x)=\chi_0(rx)$ and let $\chi_{00}$ be
the limit of $\chi_{r_k}$ of the above sequence $k\to\infty$. By Lemma \ref{chi_limit},
$\chi_0\equiv1$ in $B_r(-r, 0)$, hence
$\chi_{00}\equiv 1$ in the left half plane. Moreover we know
that $\chi_0(x)=1$ if $v_0(x)>0$, hence
$\chi_{00}\equiv1$ in the right half plane. Thus
$\chi_{00}(x)\equiv1$ in $\mathbb{R}^2$.
\par
From
$$
0=\int_{\mathbb{R}^2}|\nabla v_{00}|^2\div\phi-2\nabla
v_{00}D\phi\nabla v_{00}+\chi_{00}\div\phi\ dx
$$
for every $\phi\in C_0^\infty(\mathbb{R}^2; \mathbb{R}^2)$, we infer
that
$$
\int_{\mathbb{R}^2}|\nabla v_{00}|^2\div\phi-2\nabla
v_{00}D\phi\nabla v_{00}=0
$$
for every $\phi\in C_0^\infty(\mathbb{R}^2; \mathbb{R}^2)$. This however contradicts 
$v_{00}=c\max(0, x\cdot\nu)$. We conclude that the blow-up limit $v_{00}$
must be the constant function 0 in $\mathbb{R}^2$.\\

\par
Case 2: $0\in\partial\{v_0>0\}$.
\par
Suppose towards a contradiction that a blow-up limit
$$v_{00}=\lim_{j\to\infty} \frac{v_0(s_j\cdot)}{s_j}\not\equiv 0.$$ Similar to case 1, $v_{00}$ must be of the
form $cx\cdot\nu$ in some half plane; here $\nu$ is again a unit vector. On the other hand
Proposition \ref{conv} 5) --- which incidentally is preserved under
the blow-up limit --- implies that there are at least two vectors
$e^1\ne e^2 \in \partial B_1$ and two
half lines $\{\alpha e^1: \alpha >0\}$ and $\{\alpha e^2: \alpha >0\}$, such that
$e^1\cdot \nu>0$, $\nabla v_{00} = \nabla v_{00}\cdot e^1 \> e^1$,
$e^2\cdot \nu>0$ and $\nabla v_{00} = \nabla v_{00}\cdot e^2 \> e^2$
which is not possible for $v_{00}(x)=cx\cdot\nu$ unless $c=0$.
\end{proof}
\begin{remarks}\label{fin}
Using methods in the forthcoming paper \cite{var},
it is actually possible to show that Case 1 in the proof of Theorem \ref{main}
cannot occur. More precisely,
a degenerate point $x^0$ at which the set $\{ v_0=0\}$ contains
a disk touching $x^0$ is not possible.

\par
{\em Consequently the origin must in our example be a degenerate point.}
The fact that we can introduce many symmetry lines suggests
that we can construct degenerate points 
with growth $\le |x|^m$ in $B_{1/2}$ for any integer $m$.  
\end{remarks}
\bibliographystyle{plain}
\bibliography{zhang2.bib}
\end{document}

%% file: zhang1.pstex_t
\begin{picture}(0,0)%
\includegraphics{zhang1.pstex}%
\end{picture}%
\setlength{\unitlength}{3947sp}%
\begingroup\makeatletter\ifx\SetFigFont\undefined%
\gdef\SetFigFont#1#2#3#4#5{%
  \reset@font\fontsize{#1}{#2pt}%
  \fontfamily{#3}\fontseries{#4}\fontshape{#5}%
  \selectfont}%
\fi\endgroup%
\begin{picture}(1866,1866)(268,-1294)
\put(301,-361){\makebox(0,0)[lb]{\smash{{\SetFigFont{12}{14.4}{\rmdefault}{\mddefault}{\updefault}{\color[rgb]{0,0,0}$u>0$}%
}}}}
\put(901,239){\makebox(0,0)[lb]{\smash{{\SetFigFont{12}{14.4}{\rmdefault}{\mddefault}{\updefault}{\color[rgb]{0,0,0}$u>0$}%
}}}}
\put(901,-1036){\makebox(0,0)[lb]{\smash{{\SetFigFont{12}{14.4}{\rmdefault}{\mddefault}{\updefault}{\color[rgb]{0,0,0}$u>0$}%
}}}}
\put(1651,-361){\makebox(0,0)[lb]{\smash{{\SetFigFont{12}{14.4}{\rmdefault}{\mddefault}{\updefault}{\color[rgb]{0,0,0}$u>0$}%
}}}}
\end{picture}%

%% file: zhang2.pstex_t
\begin{picture}(0,0)%
\includegraphics{zhang2.pstex}%
\end{picture}%
\setlength{\unitlength}{3947sp}%
\begingroup\makeatletter\ifx\SetFigFont\undefined%
\gdef\SetFigFont#1#2#3#4#5{%
  \reset@font\fontsize{#1}{#2pt}%
  \fontfamily{#3}\fontseries{#4}\fontshape{#5}%
  \selectfont}%
\fi\endgroup%
\begin{picture}(1866,2016)(268,-1294)
\put(301,-361){\makebox(0,0)[lb]{\smash{{\SetFigFont{12}{14.4}{\rmdefault}{\mddefault}{\updefault}{\color[rgb]{0,0,0}$u>0$}%
}}}}
\put(901,-1036){\makebox(0,0)[lb]{\smash{{\SetFigFont{12}{14.4}{\rmdefault}{\mddefault}{\updefault}{\color[rgb]{0,0,0}$u>0$}%
}}}}
\put(1651,-361){\makebox(0,0)[lb]{\smash{{\SetFigFont{12}{14.4}{\rmdefault}{\mddefault}{\updefault}{\color[rgb]{0,0,0}$u>0$}%
}}}}
\end{picture}%

%% file: zhang3.pstex_t
\begin{picture}(0,0)%
\includegraphics{zhang3.pstex}%
\end{picture}%
\setlength{\unitlength}{3947sp}%
\begingroup\makeatletter\ifx\SetFigFont\undefined%
\gdef\SetFigFont#1#2#3#4#5{%
  \reset@font\fontsize{#1}{#2pt}%
  \fontfamily{#3}\fontseries{#4}\fontshape{#5}%
  \selectfont}%
\fi\endgroup%
\begin{picture}(1952,1952)(225,-1337)
\put(980,202){\makebox(0,0)[lb]{\smash{{\SetFigFont{9}{10.8}{\rmdefault}{\mddefault}{\updefault}{\color[rgb]{0,0,0}$u_\epsilon>\epsilon$}%
}}}}
\put(371,-404){\makebox(0,0)[lb]{\smash{{\SetFigFont{9}{10.8}{\rmdefault}{\mddefault}{\updefault}{\color[rgb]{0,0,0}$u_\epsilon>\epsilon$}%
}}}}
\put(1600,-400){\makebox(0,0)[lb]{\smash{{\SetFigFont{9}{10.8}{\rmdefault}{\mddefault}{\updefault}{\color[rgb]{0,0,0}$u_\epsilon>\epsilon$}%
}}}}
\put(992,-1044){\makebox(0,0)[lb]{\smash{{\SetFigFont{9}{10.8}{\rmdefault}{\mddefault}{\updefault}{\color[rgb]{0,0,0}$u_\epsilon>\epsilon$}%
}}}}
\end{picture}%

%% file: zhang4.pstex_t
\begin{picture}(0,0)%
\includegraphics{zhang4.pstex}%
\end{picture}%
\setlength{\unitlength}{3947sp}%
\begingroup\makeatletter\ifx\SetFigFont\undefined%
\gdef\SetFigFont#1#2#3#4#5{%
  \reset@font\fontsize{#1}{#2pt}%
  \fontfamily{#3}\fontseries{#4}\fontshape{#5}%
  \selectfont}%
\fi\endgroup%
\begin{picture}(1828,1828)(288,-1125)
\put(1435,-467){\makebox(0,0)[lb]{\smash{{\SetFigFont{10}{12.0}{\rmdefault}{\mddefault}{\updefault}{\color[rgb]{0,0,0}$x^i_t$}%
}}}}
\put(729,-13){\makebox(0,0)[lb]{\smash{{\SetFigFont{10}{12.0}{\rmdefault}{\mddefault}{\updefault}{\color[rgb]{0,0,0}$P_i$}%
}}}}
\end{picture}%

%% file: zhang5.pstex_t
\begin{picture}(0,0)%
\includegraphics{zhang5.pstex}%
\end{picture}%
\setlength{\unitlength}{3947sp}%
\begingroup\makeatletter\ifx\SetFigFont\undefined%
\gdef\SetFigFont#1#2#3#4#5{%
  \reset@font\fontsize{#1}{#2pt}%
  \fontfamily{#3}\fontseries{#4}\fontshape{#5}%
  \selectfont}%
\fi\endgroup%
\begin{picture}(1670,1154)(283,-681)
\put(1027,-626){\makebox(0,0)[lb]{\smash{{\SetFigFont{10}{12.0}{\rmdefault}{\mddefault}{\updefault}{\color[rgb]{0,0,0}$x^0$}%
}}}}
\put(1424,-240){\makebox(0,0)[lb]{\smash{{\SetFigFont{10}{12.0}{\rmdefault}{\mddefault}{\updefault}{\color[rgb]{0,0,0}$v_0\le 0$}%
}}}}
\end{picture}%

%% file: zhang2.bbl
\begin{thebibliography}{10}

\bibitem{AC}
H.~W. Alt and L.~A. Caffarelli.
\newblock Existence and regularity for a minimum problem with free boundary.
\newblock {\em J. Reine Angew. Math.}, 325:105--144, 1981.

\bibitem{AW}
J.~Andersson and G.~S. Weiss.
\newblock Cross-shaped and degenerate singularities in an unstable elliptic
  free boundary problem.
\newblock {\em J. Differential Equations}, 228(2):633--640, 2006.

\bibitem{Berestycki}
H.~Berestycki, L.~A. Caffarelli, and L.~Nirenberg.
\newblock Uniform estimates for regularization of free boundary problems.
\newblock In {\em Analysis and partial differential equations}, volume 122 of
  {\em Lecture Notes in Pure and Appl. Math.}, pages 567--619. Dekker, New
  York, 1990.

\bibitem{Bucmaster}
J.~D. Buckmaster and G.~S.~S. Ludford.
\newblock {\em Theory of laminar flames}.
\newblock Cambridge Monographs on Mechanics and Applied Mathematics. Cambridge
  University Press, Cambridge, 1982.
\newblock Electronic \& Electrical Engineering Research Studies: Pattern
  Recognition \& Image Processing Series, 2.

\bibitem{Caffarelli}
Luis~A. Caffarelli.
\newblock Uniform {L}ipschitz regularity of a singular perturbation problem.
\newblock {\em Differential Integral Equations}, 8(7):1585--1590, 1995.

\bibitem{DPS}
D.~Danielli, A.~Petrosyan, and H.~Shahgholian.
\newblock A singular perturbation problem for the {$p$}-{L}aplace operator.
\newblock {\em Indiana Univ. Math. J.}, 52(2):457--476, 2003.

\bibitem{Frehse}
J.~Frehse.
\newblock Capacity methods in the theory of partial differential equations.
\newblock {\em Jahresber. Deutsch. Math.-Verein.}, 84(1):1--44, 1982.

\bibitem{Gilbarg}
David Gilbarg and Neil~S. Trudinger.
\newblock {\em Elliptic partial differential equations of second order}, volume
  224 of {\em Grundlehren der Mathematischen Wissenschaften [Fundamental
  Principles of Mathematical Sciences]}.
\newblock Springer-Verlag, Berlin, second edition, 1983.

\bibitem{Gurevich}
Alex Gurevich.
\newblock Boundary regularity for free boundary problems.
\newblock {\em Comm. Pure Appl. Math.}, 52(3):363--403, 1999.

\bibitem{Lederman}
Claudia Lederman and Noemi Wolanski.
\newblock Viscosity solutions and regularity of the free boundary for the limit
  of an elliptic two phase singular perturbation problem.
\newblock {\em Ann. Scuola Norm. Sup. Pisa Cl. Sci. (4)}, 27(2):253--288
  (1999), 1998.

\bibitem{var}
E.~Varvaruca and G.~S. Weiss.
\newblock A geometric proof of a generalized {S}tokes conjecture.
\newblock {\em In preparation}.

\bibitem{Weiss2}
G.~S. Weiss.
\newblock A singular limit arising in combustion theory: fine properties of the
  free boundary.
\newblock {\em Calc. Var. Partial Differential Equations}, 17(3):311--340,
  2003.

\bibitem{Weiss1}
Georg~S. Weiss.
\newblock Partial regularity for weak solutions of an elliptic free boundary
  problem.
\newblock {\em Comm. Partial Differential Equations}, 23(3-4):439--455, 1998.

\end{thebibliography}
